ANNALES
DE L'INSTITUT
HENRI
POINCARÉ
PROBABILITÉS
ET STATISTIQUES

# Vitesse de convergence dans le théorème limite central pour des chaînes de Markov fortement ergodiques

## Loïc Hervé

*I.R.M.A.R., UMR-CNRS 6625, Institut National des Sciences Appliquées de Rennes, 20, Avenue des Buttes de Coüesmes CS 14 315, 35043 Rennes Cedex, France. E-mail : Loic.Herve@insa-rennes.fr.*

**Résumé.** Soit $Q$ une probabilité de transition sur un espace mesurable $E$, admettant une probabilité invariante, soit $(X_n)_n$ une chaîne de Markov associée à $Q$, et soit $\xi$ une fonction réelle mesurable sur $E$, et $S_n = \sum_{k=1}^n \xi(X_k)$. Sous des hypothèses fonctionnelles sur l'action de $Q$ et des noyaux de Fourier $Q(t)$, nous étudions la vitesse de convergence dans le théorème limite central pour la suite $(\frac{S_n}{\sqrt{n}})_n$. Selon les hypothèses nous obtenons une vitesse en $n^{-\tau/2}$ pour tout $\tau < 1$, ou bien en $n^{-1/2}$. Nous appliquons la méthode de Nagaev en l'améliorant, d'une part grâce à un théorème de perturbations de Keller et Liverani, d'autre part grâce à une majoration de $|\mathbb{E}[e^{itS_n/\sqrt{n}}] - e^{-t^2/2}|$ obtenue par une méthode de réduction en différence de martingale. Lorsque $E$ est non compact ou $\xi$ est non bornée, les conditions requises ici sur $Q(t)$ (en substance, des conditions de moment sur $\xi$) sont plus faibles que celles habituellement imposées lorsqu'on utilise le théorème de perturbation standard. Par exemple, dans le cadre des chaînes $V$-géométriquement ergodiques ou des modèles itératifs Lipschitziens, on obtient dans le t.l.c. une vitesse en $n^{-1/2}$ sous une hypothèse de moment d'ordre 3 sur $\xi$.

**Abstract.** Let $Q$ be a transition probability on a measurable space $E$ which admits an invariant probability measure, let $(X_n)_n$ be a Markov chain associated to $Q$, and let $\xi$ be a real-valued measurable function on $E$, and $S_n = \sum_{k=1}^n \xi(X_k)$. Under functional hypotheses on the action of $Q$ and the Fourier kernels $Q(t)$, we investigate the rate of convergence in the central limit theorem for the sequence $(\frac{S_n}{\sqrt{n}})_n$. According to the hypotheses, we prove that the rate is, either $\mathrm{O}(n^{-\tau/2})$ for all $\tau < 1$, or $\mathrm{O}(n^{-1/2})$. We apply the spectral Nagaev's method which is improved by using a perturbation theorem of Keller and Liverani, and a majoration of $|\mathbb{E}[e^{itS_n/\sqrt{n}}] - e^{-t^2/2}|$ obtained by a method of martingale difference reduction. When $E$ is not compact or $\xi$ is not bounded, the conditions required here on $Q(t)$ (in substance, some moment conditions on $\xi$) are weaker than the ones usually imposed when the standard perturbation theorem is used in the spectral method. For example, in the case of $V$-geometric ergodic chains or Lipschitz iterative models, the rate of convergence in the c.l.t. is $\mathrm{O}(n^{-1/2})$ under a third moment condition on $\xi$.

## 1. Introduction

Dans ce papier on désigne par $(E, \mathcal{E})$ un espace mesurable, par $Q$ une probabilité de transition sur $(E, \mathcal{E})$ admettant une probabilité invariante, notée $\nu$, par $(X_n)_{n \geq 0}$ une chaîne de Markov sur $(E, \mathcal{E})$ associée à $Q$, et enfin par $\xi$ une fonction $\nu$-intégrable de $E$ dans $\mathbb{R}$ telle que $\nu(\xi) = 0$. On pose $S_n = \sum_{k=1}^{n} \xi(X_k)$.

L'objet de ce travail est d'étudier la vitesse de convergence dans le théorème limite central (t.l.c.) pour la suite de variables aléatoires $(\xi(X_n))_{n \geq 0}$. La méthode que nous utilisons s'inspire des techniques de perturbations d'opérateurs qui ont été introduites par Nagaev [23, 24] et largement appliquées depuis, voir [10]. On trouvera dans [13] un exposé général de cette méthode et de nombreuses références.

Les hypothèses porteront sur le noyau $Q$ et les noyaux de Fourier $Q(t)$ associés à $Q$ et $\xi$ ; en substance, on supposera que, sur un certain espace de Banach, $Q$ vérifie une hypothèse de forte ergodicité et que $Q(t)$ satisfait aux conditions du théorème de perturbations d'opérateurs de Keller–Liverani [19]. Comme il apparaît déjà dans [14] (dans le cadre des modèles itératifs) et dans [15] (en vu d'obtenir un théorème local), le théorème de Keller–Liverani remplace avantageusement les énoncés classiques de perturbations d'opérateurs, notamment lorsque l'espace d'états $E$ est non compact ou que $\xi$ est non bornée. Ainsi, dans [14] et [15], nous avons pu établir des théorèmes limites sous des hypothèses de moments polynomiaux là où habituellement étaient requis des moments exponentiels.

Les résultats de vitesse de convergence dans le t.l.c. sont présentés au paragraphe 2. On appliquera ensuite (Section 3) ces résultats dans le cadre des chaînes $V$-géométriquement ergodiques et des modèles itératifs Lipschitziens.

Pour illustrer les résultats obtenus, considérons l'exemple classique sur $\mathbb{R}^d$ des processus autorégressifs $X_n = A_n X_{n-1} + b_n$ sous les hypothèses suivantes : $A_1$ est presque-sûrement une contraction stricte et $\mathbb{E}[\|b_1\|^3] < +\infty$. Alors, si $X_0$ a un moment d'ordre 2 et si $\xi$ est uniformément lipschitzienne sur $\mathbb{R}^d$ (par exemple si $\xi(x) = \|x\|$), la vitesse de convergence dans le t.l.c. est en $n^{-1/2}$ (cf. corollaire 3.2′). La condition de moment d'ordre 3 sur $b_1$ est la condition attendue en référence au théorème de Berry–Esseen pour v.a.i.i.d.

**Notations.** *Si $X$ et $Y$ sont des espaces de Banach, on désigne par $X'$ le dual topologique de $X$, par $\mathcal{L}(X)$ l'espace des endomorphismes continus de $X$, et par $\mathcal{L}(X, Y)$ l'espace des applications linéaires continues de $X$ dans $Y$. Ces espaces sont munis des normes subordonnées. On note $\langle \cdot, \cdot \rangle$ le crochet de dualité sur $X' \times X$.*

*Pour $p \geq 1$ on note $\mathbb{L}^p(\nu)$ l'espace de Lebesgue usuel associé à $\nu$. On note $\mathbf{1} = 1_E$ la fonction identiquement égale à 1 sur $E$.*

*La probabilité initiale de la chaîne sera appelée $\mu_0$. La loi normale centrée, de variance $\sigma^2$, est notée $\mathcal{N}(0, \sigma^2)$. Enfin les noyaux de Fourier associés à $Q$ et $\xi$ sont définis par*

$$t \in \mathbb{R}, x \in E, \quad Q(t)(x, dy) = e^{it\xi(y)} Q(x, dy).$$

## 2. Hypothèses et énoncés des résultats

Rappelons que $\nu$ désigne une probabilité $Q$-invariante. Les espaces de Banach sur lesquels on fera opérer $Q$ sont composés de fonctions mesurables de $E$ dans $\mathbb{C}$. Étant donné un tel espace $(\mathcal{B}, \|\cdot\|)$, nous dirons que $Q$ est $\mathcal{B}$-géométriquement ergodique si :

$\mathbf{1} \in \mathcal{B}$, $\mathcal{B}$ s'envoie continûment dans $\mathbb{L}^1(\nu)$, $Q \in \mathcal{L}(\mathcal{B})$, et il existe $\kappa_0 < 1$, $C \geq 0$ tels que

$$\forall n \geq 1, \forall f \in \mathcal{B}, \quad \|Q^n f - \nu(f)\mathbf{1}\| \leq C\kappa_0^n \|f\|.$$

Sous cette condition, si $\xi \in \mathcal{B}$, alors la série $\sum_{n \geq 0} Q^n \xi$ converge absolument dans $\mathcal{B}$ car $\nu(\xi) = 0$. Dans ce cas on pose $\check{\xi} = \sum_{n=0}^{+\infty} Q^n \xi$, et si en outre $\mathcal{B} \subset \mathbb{L}^2(\nu)$, on note

$$\sigma^2 = \nu(\check{\xi}^2) - \nu((Q\check{\xi})^2) \quad \text{et} \quad \psi = Q(\check{\xi}^2) - (Q\check{\xi})^2 - \sigma^2 \mathbf{1}.$$

**Hypothèse $(\mathcal{H})$.** *Il existe un espace de Banach $(\mathcal{B}, \|\cdot\|)$ tel que $\mathcal{B} \subset \mathbb{L}^3(\nu)$, $\xi \in \mathcal{B}$, et vérifiant en outre les conditions suivantes :*

(H1)  $Q$ *est* $\mathcal{B}$-*géométriquement ergodique.*

(H2)  *On a* $\sum_{p=0}^{+\infty} \nu(|Q^p\psi|^{3/2})^{2/3} < +\infty.$

(H3)  *On a* $\sup\{\nu(|e^{it\xi} - 1||f|), f \in \mathcal{B}, \|f\| \le 1\} = \mathrm{O}(|t|)$ *quand* $t \to 0$.

(H4)  *Il existe un intervalle ouvert* $I$ *contenant* $t = 0$ *tel que, pour* $t \in I$, *on ait* $Q(t) \in \mathcal{L}(\mathcal{B})$, *et il existe des constantes* $\kappa < 1$, $C \ge 0$ *tels que*

$$\forall n \ge 1, \forall t \in I, \forall f \in \mathcal{B}, \quad \|Q(t)^n f\| \le C\kappa^n \|f\| + C\nu(|f|),$$

*et enfin le rayon spectral essentiel de* $Q(t)$ *opérant sur* $\mathcal{B}$ *est* $\le$ *au réel* $\kappa$.

On observera que la seule condition de moment sur $\xi$ est $\nu(|\xi|^3) < +\infty$.

**Théorème 2.1.** *Supposons que* $\mu_0 = \nu$ *et* $\sigma^2 > 0$. *Sous l'hypothèse* $(\mathcal{H})$, *la suite* $(\frac{S_n}{\sqrt{n}})_n$ *converge en loi vers* $\mathcal{N}(0, \sigma^2)$ *avec une vitesse de convergence au moins en* $n^{-\tau/2}$ *pour tout* $\tau < 1$,

$$\text{à savoir : } \forall \tau < 1, \quad \sup_{x \in \mathbb{R}} \left| \mathbb{P}\left( \frac{S_n}{\sigma\sqrt{n}} \le x \right) - \mathcal{N}(0,1)(]-\infty, x]) \right| = \mathrm{O}(n^{-\tau/2}).$$

**Hypothèse** $(\widetilde{\mathcal{H}})$. *On a* $(\mathcal{H})$ *sur* $\mathcal{B}$, *avec* $\mathcal{B}$ *contenu dans un espace de Banach* $(\widetilde{\mathcal{B}}, \|\cdot\|_\sim)$ *tel que* $Q$ *est* $\widetilde{\mathcal{B}}$-*géométriquement ergodique et* $\sup\{\|Q(t)f - Qf\|_\sim, f \in \mathcal{B}, \|f\| \le 1\} = \mathrm{O}(|t|)$ *quand* $t \to 0$.

**Théorème 2.2.** *Supposons que l'hypothèse* $(\widetilde{\mathcal{H}})$ *soit satisfaite, et que* $\sigma^2 > 0$ *et* $\mu_0 \in \mathcal{B}' \cap \widetilde{\mathcal{B}}'$. *Alors la vitesse de convergence dans le t.l.c. est en* $n^{-1/2}$.

Sous la condition (H1), si $Q(\cdot)$ est de classe $\mathcal{C}^3$ de $I$ dans $\mathcal{L}(\mathcal{B})$, le théorème de perturbations standard permet d'établir une vitesse en $n^{-1/2}$ dans le t.l.c. [10]. Cependant, l'hypothèse précédente requiert que le noyau $\xi(y)^3 Q(x, dy)$ opère continûment sur $\mathcal{B}$ [10, 13], et si $\xi$ n'est pas bornée, cette condition peut être assez restrictive. En fait, ce problème apparaît déjà si l'on impose à $Q(\cdot)$ d'être continue de $I$ dans $\mathcal{L}(\mathcal{B})$. Les conditions (H3), (H4), mieux adaptées au cas $\xi$ non bornée, permettront d'appliquer le théorème de perturbations de [19].

***Remarques.***

1. *La condition* (H1) *est équivalente au fait que* $Q$ *est quasi-compact sur* $\mathcal{B}$, *avec* 1 *comme valeur propre simple et dominante.*

2. *Supposons que* $Q$ *soit* $\mathcal{B}$-*géométriquement ergodique et que* $\mathcal{B} \subset \mathbb{L}^3(\nu)$, $\xi \in \mathcal{B}$. *Soit* $(\mathcal{B}_2, \|\cdot\|_2)$ *un espace de Banach contenant les fonctions* $g^2$, $g \in \mathcal{B}$. *S'il existe* $A > 0$ *tel que l'on ait, pour* $f \in \mathcal{B}_2$, $\nu(|f|^{3/2})^{2/3} \le A\|f\|_2$, *et si* $Q$ *est* $\mathcal{B}_2$-*géométriquement ergodique, alors on a* (H2). *En effet* $\check{\xi}$, $Q\check{\xi}$ *et* **1** *sont dans* $\mathcal{B}$, *donc* $\psi \in \mathcal{B}_2$. *Par définition du nombre* $\sigma^2$ *et par invariance de* $\nu$, *on a* $\nu(\psi) = 0$, *par conséquent* $\sum_{n \ge 0} \|Q^n\psi\|_2 < +\infty$, *et* (H2) *en découle.*

3. *Soit* $(\mathcal{B}_\gamma)_{0 < \gamma \le \gamma_0}$ *une famille d'espaces de Banach, croissante pour l'inclusion. Les théorèmes 1, 2 sont particulièrement bien adaptés lorsque, pour tout* $\gamma \in ]0, \gamma_0]$, $Q$ *est* $\mathcal{B}_\gamma$-*géométriquement ergodique et que* $Q(t)$ *vérifie* (H3), (H4) *sur* $\mathcal{B}_\gamma$ *pour* $|t|$ *petit. Dans ce cas la principale difficulté réside dans le choix de paramètres* $\gamma_1 < \gamma_2 < \gamma_3$ *tels que, si* $\mathcal{B} = \mathcal{B}_{\gamma_1}$, *alors* $Q$ *vérifie sur* $\mathcal{B}_2 = \mathcal{B}_{\gamma_2}$ *les conditions de la remarque* 2 *(afin d'obtenir* (H2)*), et vérifie* $(\widetilde{\mathcal{H}})$ *avec* $\widetilde{\mathcal{B}} = \mathcal{B}_{\gamma_3}$. *Les exemples du paragraphe* 3 *font intervenir de telles familles d'espaces.*

4. *Si, pour* $t \in I$, *l'inégalité de* (H4) *est satisfaite et si l'ensemble* $Q(t)(\{\|f\| \le 1, f \in \mathcal{B}\})$ *est relativement compact dans* $(\mathcal{B}, \nu(|\cdot|))$, *alors* $Q(t)$ *est quasi-compact à itérés bornés sur* $\mathcal{B}$ [16], *et la propriété dans* (H4) *sur le rayon spectral essentiel de* $Q(t)$ *est alors automatiquement satisfaite* [11]. *Cette remarque appliquée avec* $t = 0$ *permet dans certain cas d'établir* (H1) *lorsque* $Q$ *vérifie en outre des hypothèses d'irréductibilité et d'apériodicité garantissant que* 1 *est une valeur propre simple et l'unique valeur propre de module* 1 *de* $Q$.

5. *Si $\widetilde{\mathcal{B}} = \mathcal{B}$, alors l'hypothèse $(\widetilde{\mathcal{H}})$ se réduit à $(\mathcal{H})$ et $\|Q(t) - Q\|_{\mathcal{L}(\mathcal{B})} = \mathrm{O}(|t|)$. Comme déjà mentionné, cette dernière condition est assez restrictive lorsque $\xi$ est non bornée, voir la section 3.*

6. *Une vitesse de convergence en $n^{-\tau/4}$ pour tout $\tau < 1$ a été établie dans [4] pour des chaînes de Markov stationnaires, à espace d'états compact, associées à un opérateur de transfert markovien. Dans [4] il n'est pas supposé que $Q$ a une action quasi-compacte, et les conditions sur $\xi$ sont assez faibles. Voir également [17]. Mentionnons que la fonction $\psi$ introduite au début du paragraphe est déjà utilisée dans les arguments de [4].*

7. *Grâce à (H1), la condition $\xi \in \mathcal{B}$ est commode pour définir $\check{\xi}$, et la condition $\mathcal{B} \subset \mathbb{L}^3(\nu)$ assure alors que $\check{\xi} \in \mathbb{L}^3(\nu)$. Cette dernière propriété n'est utilisée que dans la section 4.1.*

    *Si $\xi \notin \mathcal{B}$, les théorèmes 2.1, 2.2 restent vrais, à condition de renforcer la condition (H2) comme suit : il existe $\check{\xi} \in \mathbb{L}^3(\nu)$ telle que $\check{\xi} - Q\check{\xi} = \xi$, et la fonction $\psi$, que l'on peut alors définir comme au début du paragraphe, vérifie $\sum_{p=0}^{+\infty} \nu(|Q^p \psi|^{3/2})^{2/3} < +\infty$.*

La preuve des théorèmes 2.1 et 2.2, présentée en section 4, est proche de celle utilisée dans [15]. Plus précisément, on appliquera dans un premier temps, en s'inspirant de [17], une méthode de réduction en différence de martingale pour établir, sous les conditions (H1), (H2), que $\sup_{|t| \leq \sqrt{n}} |t|^{-1} |\mathbb{E}[\mathrm{e}^{\mathrm{i}tS_n/(\sigma\sqrt{n})}] - \mathrm{e}^{-t^2/2}| = \mathrm{O}(\frac{1}{\sqrt{n}})$. Grâce à cette majoration, on verra que la valeur propre dominante perturbée de $Q(t)$, fournie par (H1), (H3), (H4) et les résultats de [19], vérifie le développement limité requis pour l'application des techniques usuelles de transformée de Fourier.

## 3. Exemples

Dans les deux exemples traités dans ce paragraphe, on désigne par $(E, d)$ un espace métrique non compact tel que toute boule fermée de $E$ soit compacte. On munit $E$ de sa tribu borélienne $\mathcal{E}$, et l'on note $x_0$ un point quelconque de $E$.

### 3.1. Application aux chaînes $V$-géométriquement ergodiques

Soit $V$ une fonction mesurable de $E$ dans $[1, +\infty[$ telle que $V(x) \to +\infty$ quand $d(x, x_0) \to +\infty$, et soit $(X_n)_{n \geq 0}$ une chaîne $V$-géométriquement ergodique [21] (chap. 16), à savoir : il existe une probabilité $Q$-invariante, $\nu$, telle que $\nu(V) < +\infty$, et $Q$ est $\mathcal{B}_V$-géométriquement ergodique, où $\mathcal{B}_V$ est l'espace des fonctions mesurables sur $E$, à valeurs complexes, vérifiant

$$\|f\|_V = \sup\{V(x)^{-1}|f(x)|, x \in E\} < +\infty.$$

On supposera en outre que la tribu $\mathcal{E}$ est à base dénombrable, et que $(X_n)_{n \geq 0}$ est apériodique et $\psi$-irréductible vis-à-vis d'une certaine mesure positive $\psi$ $\sigma$-finie. Rappelons que, si $\xi^2$ est dominée par un multiple de $V$, alors $(\frac{S_n}{\sqrt{n}})_n$ converge en loi vers $\mathcal{N}(0, \sigma^2)$ [21].

**Corollaire 3.1.** *Si $|\xi|^3$ est dominée par un multiple de $V$, si $\sigma^2 > 0$ et si $\mathbb{E}[V(X_0)] < +\infty$, alors $(\xi(X_n))_{n \geq 0}$ vérifie le t.l.c. avec une vitesse en $n^{-1/2}$.*

On observera que ce corollaire s'applique en particulier quand la loi initiale est la masse de Dirac $\delta_x$ en un point quelconque $x$ de $E$.

Les chaînes géométriquement ergodiques étant fortement mélangeantes [21] (chap. 16), voir également [18], la conclusion du corollaire peut être établie, dans le cas stationnaire, grâce au théorème de Bolthausen [3] sous l'hypothèse $\nu(|\xi|^{3+\varepsilon}) < +\infty$, et l'on peut également appliquer [6] qui étudie la vitesse de convergence en norme $\mathbb{L}^1$ dans le t.l.c.

Dans le cadre spécifique des chaînes $V$-géométriquement ergodiques, on trouvera dans [8], sous des conditions de moment assez fortes, des résultats de vitesse de convergence dans le t.l.c. (cas stationnaire) sous forme d'inégalités de Paley, et dans [20] des développements d'Edgeworth sous la condition $\xi$ bornée. Enfin [25] établit une vitesse en $(\frac{\ln n}{n})^\beta$ lorsque $\xi$ est dominée par un multiple de $V^\alpha$ $(0 < \alpha \leq \frac{1}{2})$, avec $\beta = \frac{1}{2(\alpha+1)}$.

**Démonstration du corollaire 3.1.** Soit $W = V^{1/3}$, $U = V^{2/3}$, et soient $\mathcal{B}_W$, $\mathcal{B}_U$ les espaces obtenus comme ci-dessus en remplaçant $V$ respectivement par $W$ et $U$. Notons que $\nu \in \mathcal{B}'_W$ car $\nu(W) < +\infty$, que $\xi \in \mathcal{B}_W$, et enfin que $\mathcal{B}_W$ s'envoie continûment dans $\mathbb{L}^3(\nu)$ car $\nu(W^3) = \nu(V) < +\infty$. Nous allons montrer que $Q$ vérifie l'hypothèse $(\widetilde{\mathcal{H}})$ avec $\mathcal{B} = \mathcal{B}_W$, $\widetilde{\mathcal{B}} = \mathcal{B}_V$.

Comme $Q$ est $V$-géométriquement ergodique par hypothèse, il l'est également relativement à $W$, voir [21] (lemma 15.2.9 et theorem 16.0.1). D'où (H1). Pour les mêmes raisons $Q$ est $U$-géométriquement ergodique. En outre $\mathcal{B}_U$ s'envoie continûment dans $L^{3/2}(\nu)$ car $\nu(U^{3/2}) = \nu(V) < +\infty$. D'où (H2) (cf. remarque 2 avec $\mathcal{B}_2 = \mathcal{B}_U$). Par ailleurs (H3) résulte de l'inégalité $|e^{it\xi} - 1| \leq |t||\xi|$, et on a clairement $Q(t) \in \mathcal{L}(\mathcal{B}_W)$ pour tout $t \in \mathbb{R}$. La propriété (H4), établie en appendice, résulte d'un travail récent de Hennion [12]. Enfin si $f \in \mathcal{B}_W$, alors

$$|Q(t)f - Qf| \leq Q(|e^{it\xi} - 1||f|) \leq |t|Q(|\xi|W)\|f\|_W \leq |t|\|\xi\|_W Q(W^2)\|f\|_W,$$

d'où $\|Q(t)f - Qf\|_V \leq |t| \|\xi\|_W \|QW^2\|_V \|f\|_W$. Ce qui précède prouve $(\widetilde{\mathcal{H}})$. Comme par hypothèse $\mu_0(V) < +\infty$, on a $\mu_0 \in \mathcal{B}'_W \cap \mathcal{B}'_V$, et le corollaire découle du théorème 2.2. $\qquad \square$

### 3.2. Application aux modèles itératifs lipschitziens

On désigne par $G$ un semi-groupe de transformations lipschitziennes de $E$, et par $\mathcal{G}$ une tribu sur $G$. On suppose que l'action de $G$ sur $E$ est mesurable. Pour $g \in G$, on pose

$$c(g) = \sup\left\{\frac{d(gx, gy)}{d(x,y)}, x, y \in E, x \neq y\right\}.$$

Soit $(Y_n)_{n\geq 1}$ une suite de v.a.i.i.d. à valeurs dans $G$. On note $\pi$ leur loi commune. Étant donnée une variable aléatoire $X_0$ à valeurs dans $E$, indépendante de la suite $(Y_n)_n$, de loi $\mu_0$, on considère la suite $(X_n)_{n\geq 0}$ définie pour tout $n \geq 1$ par $X_n = Y_n X_{n-1}$.

Alors $(X_n)_{n\geq 0}$ est une chaîne de Markov de probabilité de transition $(Qf)(x) = \int_G f(gx) \, d\pi(g)$.

On suppose dans ce paragraphe qu'il existe une constante $C \geq 0$ telle que l'on ait

$$\forall(x,y) \in E^2, \quad |\xi(x) - \xi(y)| \leq Cd(x,y).$$

Supposons que $\int_G c(g)^2 \, d\pi(g) < 1$ et $\int_G d(gx_0, x_0)^2 \, d\pi(g) < +\infty$. Alors la suite $(\frac{S_n}{\sqrt{n}})_n$ converge en loi vers une gaussienne $\mathcal{N}(0, \sigma^2)$, voir [2].

Soit $\Gamma(g) = 1 + c(g) + d(gx_0, x_0)$. On renforce ici la condition précédente en supposant qu'il existe un entier $n_0 \geq 1$ tel que

$$\int_G \Gamma(g)^3(1 + c(g)^{1/2}) \, d\pi(g) < +\infty \quad \text{et} \quad \int_G c(g)^{1/2} \max\{c(g), 1\}^3 \, d\pi^{*n_0}(g) < 1, \tag{$*$}$$

où l'on a désigné par $\pi^{*n_0}$ la loi de $Y_{n_0} \cdots Y_1$.

Sous ces hypothèses, on sait qu'il existe une unique probabilité $Q$-invariante, $\nu$, et que $\nu(d(\cdot, x_0)^3) < +\infty$ (appliquer le théorème I de [14] avec la distance $d(x,y)^{1/2}$).

**Corollaire 3.2.** *Sous l'hypothèse* $(*)$, *si* $\sigma^2 > 0$ *et* $\mathbb{E}[d(X_0, x_0)^2] < +\infty$, *alors* $(\xi(X_n))_{n\geq 0}$ *vérifie le t.l.c. avec une vitesse en* $n^{-1/2}$.

Le cadre ci-dessus contient celui des modèles itératifs [7], en particulier celui des processus autorégressifs définis par une v.a $X_0$ à valeurs dans $\mathbb{R}^d$, puis par $X_n = A_n X_{n-1} + b_n$ $(n \geq 1)$, où $(A_n, b_n)_{n\geq 1}$ est une suite de v.a.i.i.d. à valeurs dans $\mathcal{M}_d(\mathbb{R}) \times \mathbb{R}^d$, indépendante de $X_0$ ; on a noté $\mathcal{M}_d(\mathbb{R})$ l'espace des matrices réelles carrées d'ordre $d$.

Dans ce contexte une vitesse en $\frac{1}{\sqrt{n}}$ dans le t.l.c. a été établie dans [22] lorsque $A_1$ est presque-sûrement une contraction stricte et que $b_1$ vérifie une condition de moment exponentiel. Dans [5] une vitesse en $n^{-\tau/2}$ $(\tau < 1)$ est obtenue pour une large classe de fonctions $\xi$, sous la même condition sur $A_1$ et sous l'hypothèse $\mathbb{E}[\|b_1\|^p] < +\infty$ pour tout $p \in \mathbb{N}$. Les résultats de [22] ont été étendus dans [14] sous une condition de contraction en moyenne sur $\|A_1\|$ et sous la condition de moment $\mathbb{E}[\|b_1\|^{8+\varepsilon}] < +\infty$ $(\varepsilon > 0)$, où $\|\cdot\|$ désigne indifféremment une norme de $\mathbb{R}^d$ et la norme subordonnée associée sur $\mathcal{M}_d(\mathbb{R})$. Du corollaire 3.2 nous déduisons par exemple le résultat suivant.

**Corollaire 3.2′.** *Si $\|A_1\| < 1$ presque sûrement, si $\mathbb{E}[\|b_1\|^3] < +\infty$ et $\mathbb{E}[\|X_0\|^2] < +\infty$, alors $(\xi(X_n))_{n \geq 0}$ vérifie le t.l.c. avec une vitesse en $n^{-1/2}$ (sous réserve que $\sigma^2 > 0$).*

Supposons $A_n = 0$, $d = 1$, et $\xi(x) = x - \mathbb{E}[b_1]$. Alors $X_n = b_n$ pour tout $n \geq 1$, donc $S_n = (b_1 + \cdots + b_n - n\mathbb{E}[b_1])/n$, et nous retrouvons ainsi la même conclusion que le théorème de Berry–Esseen, sous la même condition de moment d'ordre 3.

**Démonstration du corollaire 3.2.** Pour simplifier nous supposons que $n_0 = 1$ dans la condition $(*)$. Soit $\lambda \in \,]0, 1]$ quelconque, soit $p_\lambda(x) = 1 + \lambda d(x, x_0)^{1/2}$ (le réel $\lambda$ sera choisi ultérieurement), et pour $\gamma > 0$, soit $\mathcal{L}_\gamma$ l'espace de Banach des fonctions $f$ de $E$ dans $\mathbb{C}$ telles que

$$m_\gamma(f) = \sup\left\{ \frac{|f(x) - f(y)|}{d(x,y)^{1/2} p_\lambda(x)^\gamma p_\lambda(y)^\gamma}, x, y \in E, x \neq y \right\} < +\infty,$$

muni de la norme $\|f\|_\gamma = m_\gamma(f) + |f|_\gamma$, où l'on a posé $|f|_\gamma = \sup_{x \in E}(|f(x)|/p_\lambda(x)^{\gamma+1})$.

On a clairement $\xi \in \mathcal{L}_1$, et si $f \in \mathcal{L}_1$, alors $\nu(|f|^3)^{1/3} \leq |f|_1 \nu(p_\lambda^6)^{1/3}$, avec $\nu(p_\lambda^6) < +\infty$. Donc $\mathcal{L}_1$ s'envoie continûment dans $\mathbb{L}^3(\nu)$.

**Lemme 3.1.** *On a* (H1) (H3) *avec* $\mathcal{B} = \mathcal{L}_1$, *et* (H2).

**Démonstration.** (H3) résulte de $|\mathrm{e}^{\mathrm{i}t\xi} - 1| \leq |t||\xi| \leq |t|[|\xi(x_0)| + Cd(x, x_0)]$ et de la définition de $|\cdot|_1$. On a $\nu \in \mathcal{L}'_1 \cap \mathcal{L}'_3$ car $\nu(p_\lambda^4) < +\infty$. En choisissant $\lambda$ suffisamment petit, l'ergodicité géométrique de $Q$ relativement à $\mathcal{B} = \mathcal{L}_1$, puis à $\mathcal{L}_3$, résulte des conditions $(*)$ et de [14] (théorème 5.5 appliqué avec la distance $d(x, y)^{1/2}$).[1] En outre si $f \in \mathcal{L}_3$, alors $\nu(|f|^{3/2})^{2/3} \leq |f|_3 \nu(p_\lambda^6)^{2/3} \leq \|f\|_3 \nu(p_\lambda^6)^{2/3}$, et $\mathcal{L}_3$ contient les fonctions $g^2$, $g \in \mathcal{L}_1$. D'où (H2) (cf. remarque 2 avec $\mathcal{B}_2 = \mathcal{L}_3$). □

**Lemme 3.2.** *$Q(t)$ est un endomorphisme continu de $\mathcal{L}_1$ pour tout $t \in \mathbb{R}$, et enfin on a* (H4) *pour $|t|$ petit.*

**Démonstration.** On pose $\delta_\lambda(g) = \max\{c(g), 1\}^{1/2} + \lambda d(gx_0, x_0)^{1/2}$. On a $\sup_{x \in E}(p_\lambda(gx)/p_\lambda(x)) \leq \delta_\lambda(g)$ et $\delta_\lambda(g)^2 \leq 2\Gamma(g)$. Soit $f \in \mathcal{L}_1$. De la définition de $c(g)$ et de l'inégalité $|\mathrm{e}^{\mathrm{i}a} - \mathrm{e}^{\mathrm{i}b}| \leq 2|b - a|^{1/2}$ $(a, b \in \mathbb{R})$, on obtient en posant $A = \pi(c^{1/2}\delta_\lambda^2)$

$$|(Q(t)f)(x) - (Q(t)f)(y)| \leq \int |f(gx) - f(gy)|\,d\pi(g) + \int |f(gy)||\mathrm{e}^{\mathrm{i}t\xi(gx)} - \mathrm{e}^{\mathrm{i}t\xi(gy)}|\,d\pi(g)$$

$$\leq A\,m_1(f)d(x,y)^{1/2}p_\lambda(x)p_\lambda(y) + 2AC^{1/2}|t|^{1/2}|f|_1 d(x,y)^{1/2}p_\lambda(x)p_\lambda(y)^2.$$

On peut évidemment supposer que $d(y, x_0) \leq d(x, x_0)$ (sinon, inverser le rôle joué par $x$ et $y$), de sorte que $p_\lambda(y)^2 \leq p_\lambda(x)p_\lambda(y)$, et il vient que $Q(t)f \in \mathcal{L}_1$, avec

$$m_1(Q(t)f) \leq Am_1(f) + 2AC^{1/2}|t|^{1/2}|f|_1.$$

---

[1] À cet effet le réel $\lambda$ doit être fixé tel que $\pi(c^{1/2}\delta_\lambda^6) < 1$, où $\delta_\lambda(g) = \max\{c(g), 1\}^{1/2} + \lambda d(gx_0, x_0)^{1/2}$, ce qui est possible grâce aux conditions $(*)$ et au théorème de convergence dominée.

On a clairement $|Q(t)f|_1 \le |f|_1 \pi(\delta_\lambda^2)$, de sorte que $Q(t)$ a une action continue sur $\mathcal{L}_1$.

Pour établir l'inégalité de (H4), on utilise le fait que les normes $\|\cdot\|_1$ et $\|\cdot\|_\nu = m_1(\cdot) + \nu(|\cdot|)$ sont équivalentes [14] (Section 5). Alors, d'après l'inégalité ci-dessus, il existe une constante $D > 0$ telle que

$$m_1(Q(t)f) \le Am_1(f) + D|t|^{1/2}[m_1(f) + \nu(|f|)] = (A + D|t|^{1/2})m_1(f) + D|t|^{1/2}\nu(|f|).$$

On a $\nu(|Q(t)f|) \le \nu(Q|f|) = \nu(|f|)$. Donc $\|Q(t)f\|_\nu \le (A + D|t|^{1/2})\|f\|_\nu + (D|t|^{1/2} + 1)\nu(|f|)$. Le réel $\lambda$ fixé dans la preuve du lemme 3.1 est tel que $A = \pi(c^{1/2}\delta_\lambda^2) < 1$. Soit alors $t_0$ tel que $\kappa = A + D|t_0|^{1/2} < 1$, soit $C' = D|t_0|^{1/2} + 1$, et soit $t$ tel que $|t| \le t_0$. Il résulte d'une récurrence évidente que $\|Q(t)^n f\|_\nu \le \kappa^n \|f\|_\nu + \frac{C'}{1-\kappa}\nu(|f|)$, ce qui prouve l'inégalité de (H4).

Il reste à établir dans (H4) la propriété relative au rayon spectral essentiel de $Q(t)$. On a $\nu(|Q(t)^n f|) \le \nu(Q^n|f|) = \nu(|f|)$, et la boule unité de $(\mathcal{L}_1, \|\cdot\|_\nu)$ est relativement compacte dans $(\mathcal{L}_1, \nu(|\cdot|))$ (utiliser le théorème d'Ascoli et le théorème de Lebesgue, cf. [14], Lem. 5.4). La propriété souhaitée résulte alors de l'inégalité de (H4) et de [11].                                                                                □

On a déjà vu que $Q$ est $\mathcal{L}_3$-géométriquement ergodique. Enfin on a :

**Lemme 3.3.** *Il existe $E > 0$ tel que, pour $f \in \mathcal{L}_1$, $t \in \mathbb{R}$, on ait $\|Q(t)f - Qf\|_3 \le E|t|\|f\|_1$.*

**Démonstration.** Soit $f \in \mathcal{L}_1$. Il existe $D > 0$ tel que $|\xi| \le D p_\lambda^2$. D'où

$$|(Q(t)f)(x) - Qf(x)| \le \int |e^{it\xi(gx)} - 1||f(gx)|\, d\pi(g) \le D|f|_1|t|p_\lambda(x)^4\pi(\delta_\lambda^4),$$

donc $|Q(t)f - Qf|_3 \le D|f|_1|t|\pi(\delta_\lambda^4)$. En outre, en posant $\tilde{A} = \pi(c^{1/2}\delta_\lambda^4)$, $\tilde{B} = \pi(c\delta_\lambda^2)$, on a

$$|[(Q(t)f)(x) - Qf(x)] - [(Q(t)f)(y) - Qf(y)]|$$
$$\le \tilde{A}Dm_1(f)|t|d(x,y)^{1/2}p_\lambda(x)^3p_\lambda(y) + \tilde{B}C|t||f|_1d(x,y)p_\lambda(y)^2.$$

On a $p_\lambda(y) \le p_\lambda(y)^3$ et $d(x,y)^{1/2} \le \frac{1}{\lambda}(p_\lambda(x) + p_\lambda(y)) \le \frac{2}{\lambda}p_\lambda(x)p_\lambda(y) \le \frac{2}{\lambda}p_\lambda(x)^3p_\lambda(y)$, par conséquent $m_3(Q(t)f - Qf) \le (\tilde{A}Dm_1(f) + \frac{2}{\lambda}\tilde{B}C|f|_1)|t|$.                                                                                □

Les lemmes précédents montrent que $Q(\cdot)$ vérifie $(\widetilde{\mathcal{H}})$ avec $\mathcal{B} = \mathcal{L}_1$ et $\widetilde{\mathcal{B}} = \mathcal{L}_3$. Enfin, comme $\mu_0(p_\lambda^4) < +\infty$ par hypothèse, on a $\mu_0 \in \mathcal{L}_1' \cap \mathcal{L}_3'$. Le corollaire 3.2 résulte du théorème 2.2.                                                                                □

# 4. Démonstration des théorèmes 2.1, 2.2

## 4.1. Une inégalité sur les fonctions caractéristiques

On suppose dans ce paragraphe que $X_0$ suit la loi $\nu$, que les conditions (H1), (H2) sont satisfaites, et enfin que $\mathcal{B} \subset \mathbb{L}^3(\nu)$, $\xi \in \mathcal{B}$. Les éléments $\check{\xi}$, $\sigma^2$ et $\psi$ ont été définis au début du paragraphe 2, et pour $n \ge 1$ on pose

$$U_n = \check{\xi}(X_n) - Q\check{\xi}(X_{n-1}) \quad \text{et} \quad T_n = U_1 + \cdots + U_n.$$

**Proposition 4.1.** *Si $\sigma^2 > 0$, alors il existe une constante $C > 0$ telle que l'on ait*

$$\forall n \in \mathbb{N}^*, \forall t \in [-\sqrt{n}, \sqrt{n}], \quad |\mathbb{E}[e^{itT_n/(\sigma\sqrt{n})}] - e^{-t^2/2}| \le C\frac{|t|}{\sqrt{n}}.$$

Admettons pour le moment cette proposition. En utilisant la méthode de Gordin [9], nous allons en déduire le résultat suivant.

**Corollaire 4.1.** *Si* $\sigma^2 > 0$, *alors il existe une constante* $C > 0$ *telle que l'on ait*

$$\forall n \in \mathbb{N}^*, \ \forall t \in [-\sqrt{n}, \sqrt{n}], \quad |\mathbb{E}[e^{itS_n/(\sigma\sqrt{n})}] - e^{-t^2/2}| \le C \frac{|t|}{\sqrt{n}}.$$

**Démonstration.** Soit $V_n = Q\check{\xi}(X_0) - Q\check{\xi}(X_n)$. Grâce à l'équation de Poisson $\check{\xi} - Q\check{\xi} = \xi$, on obtient que $S_n = T_n + V_n$. Par ailleurs, de la stationnarité de $(X_n)_n$ et du fait que $\check{\xi} \in \mathcal{B} \subset \mathbb{L}^1(\nu)$, il vient que $\sup_n \mathbb{E}[|V_n|] < +\infty$. Enfin on a

$$|\mathbb{E}[e^{itS_n/(\sigma\sqrt{n})}] - e^{-t^2/2}| = |\mathbb{E}[e^{itT_n/(\sigma\sqrt{n})} e^{itV_n/(\sigma\sqrt{n})}] - e^{-t^2/2}|$$
$$\le |\mathbb{E}[e^{itT_n/(\sigma\sqrt{n})}] - e^{-t^2/2}| + \mathbb{E}[|e^{itV_n/(\sigma\sqrt{n})} - 1|],$$

avec $\mathbb{E}[|e^{itV_n/(\sigma\sqrt{n})} - 1|] \le \frac{1}{\sigma} \frac{|t|}{\sqrt{n}} \sup_n \mathbb{E}[|V_n|]$. On conclut alors grâce à la proposition. $\qquad\square$

**Démonstration de la proposition 4.1.** On note $\mathcal{F}_n = \sigma(X_0, \dots, X_n)$ pour $n \ge 0$. Alors $(U_n)_{n\ge 1}$ est une suite stationnaire d'accroissements de martingale relativement à $(\mathcal{F}_n)_n$, et $\mathbb{E}[|U_1|^3] < +\infty$ car $\check{\xi} \in \mathcal{B} \subset \mathbb{L}^3(\nu)$.

Pour simplifier nous considérons le cas $\sigma^2 = 1$, et nous adaptons au cadre markovien les majorations présentées dans la preuve du théorème 6 de [17] (pp. 41–44). À cet effet on pose $T_0 = 0$, $W_n = U_n^2 - 1$ pour $n \ge 1$, et l'on rappelle que $e^{ix} = 1 + ix - \frac{x^2}{2} + u(ix)$, avec $|u(ix)| \le \frac{|x|^3}{6}$.

En écrivant $\mathbb{E}[e^{i\lambda T_n}] = \mathbb{E}[e^{i\lambda T_{n-1}} e^{i\lambda U_n}]$, puis en appliquant la remarque précédente avec $x = \lambda U_n$, et enfin en observant que $\mathbb{E}[e^{i\lambda T_{n-1}} U_n] = \mathbb{E}[e^{i\lambda T_{n-1}} \mathbb{E}[U_n | \mathcal{F}_{n-1}]] = 0$, il est facile de voir par récurrence que pour $t \in \mathbb{R}$ et $n \ge 1$

$$\mathbb{E}[e^{i(t/\sqrt{n})T_n}] - e^{-t^2/2} = A_n(t) + B_n(t) + C_n(t) \quad \text{avec}$$

$$A_n(t) = \left(1 - \frac{t^2}{2n}\right)^n - e^{-t^2/2}, \quad \text{puis} \quad B_n(t) = \sum_{k=0}^{n-1} \left(1 - \frac{t^2}{2n}\right)^k \mathbb{E}\left[e^{i(t/\sqrt{n})T_{n-k-1}} u\left(i\frac{t}{\sqrt{n}} U_{n-k}\right)\right],$$

$$\text{et enfin} \quad C_n(t) = -\frac{t^2}{2n} \sum_{k=0}^{n-1} \left(1 - \frac{t^2}{2n}\right)^k \mathbb{E}[e^{i(t/\sqrt{n})T_{n-k-1}} W_{n-k}].$$

Soit $t \in [-\sqrt{n}, \sqrt{n}]$. On a

$$0 \le -A_n(t) \le e^{-t^2/2} - e^{n(-t^2/(2n) - 3t^4/(8n^2))} \le (1 - e^{-3t^4/(8n)}) e^{-t^2/2} \le \frac{3t^4}{8n} e^{-t^2/2} \le C_1 \frac{|t|}{\sqrt{n}},$$

$$|B_n(t)| \le \sum_{k=0}^{n-1} \left(1 - \frac{t^2}{2n}\right)^k \frac{|t|^3}{6n\sqrt{n}} \mathbb{E}[|U_{n-k}|^3] = \frac{1}{3} \mathbb{E}[|U_1|^3] \frac{|t|}{\sqrt{n}}.$$

Pour la majoration de $C_n(t)$, on utilise le lemme suivant.

**Lemme 4.1.** *Pour* $k \ge \ell \ge 1$ *on a* $\mathbb{E}[W_k | \mathcal{F}_{\ell-1}] = (Q^{k-\ell} \psi)(X_{\ell-1})$, *où* $\psi$ *est la fonction définie au début de la section 2.*

**Démonstration.** On a $U_k^2 = \check{\xi}^2(X_k) - 2\check{\xi}(X_k) Q\check{\xi}(X_{k-1}) + (Q\check{\xi})^2(X_{k-1})$. Comme $(X_n)_{n\ge 0}$ est une chaîne de Markov, on a $\mathbb{E}[\check{\xi}(X_k) Q\check{\xi}(X_{k-1}) | \mathcal{F}_{k-1}] = \mathbb{E}[(Q\check{\xi})^2(X_{k-1}) | \mathcal{F}_{k-1}]$, puis

$$\mathbb{E}[U_k^2 | \mathcal{F}_{\ell-1}] = Q^{k-\ell+1} \check{\xi}^2(X_{\ell-1}) - Q^{k-\ell}(Q\check{\xi})^2(X_{\ell-1}).$$

D'où $\mathbb{E}[W_k^2 | \mathcal{F}_{\ell-1}] = \mathbb{E}[U_k^2 | \mathcal{F}_{\ell-1}] - 1 = Q^{k-\ell}(Q\check{\xi}^2 - (Q\check{\xi})^2 - \mathbf{1})(X_{\ell-1}) = Q^{k-\ell} \psi(X_{\ell-1})$. $\qquad\square$

D'après (H2) on sait que $C_3 = \sum_{p=0}^{+\infty} \nu(|Q^p\psi|^{3/2})^{2/3} < +\infty$, donc $\sum_{p=0}^{+\infty} \nu(|Q^p\psi|) < +\infty$. Effectuons maintenant sur $W_\ell$ une réduction en différence de martingale, à savoir $W_\ell = Y_\ell + Z_\ell$ pour $\ell \geq 1$, avec

$$Y_\ell = \sum_{p=0}^{+\infty} [\mathbb{E}[W_{p+\ell}|\mathcal{F}_\ell] - \mathbb{E}[W_{p+\ell}|\mathcal{F}_{\ell-1}]] \quad \text{et} \quad Z_\ell = \sum_{p=0}^{+\infty} \mathbb{E}[W_{p+\ell}|\mathcal{F}_{\ell-1}] - \sum_{p=1}^{+\infty} \mathbb{E}[W_{p+\ell}|\mathcal{F}_\ell].$$

$Y_\ell$ est $\mathcal{F}_\ell$-mesurable, et $\mathbb{E}[Y_\ell|\mathcal{F}_{\ell-1}] = 0$. Comme $T_{\ell-1}$ est $\mathcal{F}_{\ell-1}$-mesurable, on obtient $\mathbb{E}[e^{itT_{\ell-1}}Y_\ell] = \mathbb{E}[e^{itT_{\ell-1}} \times \mathbb{E}[Y_\ell|\mathcal{F}_{\ell-1}]] = 0$, d'où $\mathbb{E}[e^{itT_{\ell-1}}W_\ell] = \mathbb{E}[e^{itT_{\ell-1}}Z_\ell]$. Par conséquent, en posant $Z'_\ell = \sum_{p=0}^{+\infty} \mathbb{E}[W_{p+\ell}|\mathcal{F}_{\ell-1}]$, il vient

$$\mathbb{E}[e^{itT_{\ell-1}}W_\ell] = \mathbb{E}[e^{itT_{\ell-1}}Z'_\ell] - \mathbb{E}[e^{itT_{\ell-1}}Z'_{\ell+1}]. \qquad (E)$$

On a $Z'_\ell = \sum_{p=0}^{+\infty} (Q^p\psi)(X_{\ell-1})$ (lemme 4.1), et comme $\sum_{p=0}^{+\infty} \nu(|Q^p\psi|) < +\infty$, il vient

$$\mathbb{E}[e^{itT_{\ell-1}}Z'_\ell] = \sum_{p=0}^{+\infty} \mathbb{E}[e^{itT_{\ell-1}}(Q^p\psi)(X_{\ell-1})]$$

$$= \sum_{p=0}^{+\infty} \mathbb{E}[e^{it(U_2+\cdots+U_\ell)}(Q^p\psi)(X_\ell)] \quad \text{(par stationnarité de } (X_n)_{n\geq 0})$$

$$= \sum_{p=0}^{+\infty} \mathbb{E}[e^{it(T_\ell-U_1)}\mathbb{E}[W_{p+\ell+1}|\mathcal{F}_\ell]] = \mathbb{E}[e^{it(T_\ell-U_1)}Z'_{\ell+1}] \quad \text{(par le lemme 4.1).}$$

De l'égalité ($E$) on déduit que

$$|\mathbb{E}[e^{itT_{\ell-1}}W_\ell]| \leq \mathbb{E}[|Z'_{\ell+1}||e^{it(U_\ell-U_1)} - 1|] \leq |t|\mathbb{E}[|Z'_{\ell+1}|^{3/2}]^{2/3}\mathbb{E}[|U_\ell - U_1|^3]^{1/3}$$

avec $\mathbb{E}[|Z'_\ell|^{3/2}]^{2/3} \leq C_3$, puis par stationnarité $\mathbb{E}[|U_\ell - U_1|^3]^{1/3} \leq 2\mathbb{E}[|U_1|^3]^{1/3}$. Par conséquent $|\mathbb{E}[e^{it/(\sqrt{n})T_{\ell-1}}W_\ell]| \leq 2\mathbb{E}[|U_1|^3]^{1/3}C_3\frac{|t|}{\sqrt{n}}$, et finalement $|C_n(t)| \leq 2C_3\mathbb{E}[|U_1|^3]^{1/3}\frac{|t|}{\sqrt{n}}$. $\qquad \square$

### 4.2. Un théorème de perturbations

**Théorème 4.1.** *Supposons que les conditions* (H1), (H3), (H4) *soient satisfaites.*

*Soit $0 < \tau < 1$. Il existe un intervalle ouvert $J \subset I$ centré en $t = 0$, et des applications $\lambda(\cdot)$, $v(\cdot)$, $\phi(\cdot)$ et $N(\cdot)$ à valeurs respectivement dans $\mathbb{C}$, $\mathcal{B}$, $\mathcal{B}'$ et $\mathcal{L}(\mathcal{B})$ tels que l'on ait, pour $t \in J$, $n \geq 1$, et $f \in \mathcal{B}$,*

$$Q(t)^n f = \lambda(t)^n \langle \phi(t), f \rangle v(t) + N(t)^n f, \qquad (D)$$

*avec en outre les propriétés suivantes :*

(a) $\langle \phi(t), v(t) \rangle = 1$, $\phi(t)N(t) = 0$, $N(t)v(t) = 0$, $Q(t)v(t) = \lambda(t)v(t)$, *et* $\lim_{t \to 0} \lambda(t) = 1$.
(b) $\langle \nu, v(t) \rangle = 1$, *et il existe $C > 0$ et $\rho < 1$ tels que l'on ait pour $t \in J$ et $n \geq 1$*
    (b1) $\langle \nu, |v(t) - \mathbf{1}| \rangle \leq C|t|^\tau$,
    (b2) $|\langle \phi(t), \mathbf{1} \rangle - 1| \leq C|t|^\tau$,
    (b3) $\langle \nu, |N(t)^n \mathbf{1}| \rangle \leq C\rho^n |t|^\tau$.
(c) $\lambda(\cdot)$ *est continue sur $J$.*

**Démonstration.** Soit $t \in I$, $n \geq 1$, $f \in \mathcal{B}$. On a $|Q(t)^n f| \leq Q^n|f|$, donc par invariance de $\nu$, $\nu(|Q(t)^n f|) \leq \nu(Q^n|f|) = \nu(|f|)$, puis pour $t_0 \in I$ et $h \in \mathbb{R}$ tel que $t_0 + h \in I$

$$\nu(|Q(t_0 + h)f - Q(t_0)f|) \leq \nu(Q(|e^{ih\xi} - 1||f|)) = \nu(|e^{ih\xi} - 1||f|).$$

D'où $\sup\{\nu(|Q(t_0+h)f - Q(t_0)f|),\ f \in \mathcal{B},\ \|f\| \le 1\} = \mathrm{O}(|h|)$ d'après (H3).

Les deux propriétés précédentes et la condition (H4) permettent d'appliquer le théorème de Keller–Liverani [1, 19].[2] Pour les assertions (a) et (b) nous appliquons ce théorème en $t_0 = 0$. Sous la condition (H1), celui-ci assure la décomposition (D) et l'assertion (a). Le point (b) est une conséquence de résultats intermédiaires contenus dans [19] dont nous rappelons les principaux arguments.

Précisons que le réel appelé $r$ dans [19], théorème 1, est choisi ici tel que $\frac{\ln \kappa - \ln r}{\ln \kappa} = \tau$, où $\kappa$ est le réel de (H4). Donc $0 < r < 1$. Soit $\Gamma_1$ (resp. $\Gamma_0$) un cercle orienté de centre $z = 1$, de rayon suffisamment petit (resp. de centre $z = 0$, de rayon $\rho$ tel que $r < \rho < 1$). D'après [19], théorème 1, il existe $C > 0$ telle que l'on ait pour $z \in \Gamma_0 \cup \Gamma_1$, $t \in J$, $f \in \mathcal{B}$,

$$\nu(|(z-Q(t))^{-1}f| - (z-Q)^{-1}f|) \le C|t|^{\tau}\|f\|. \qquad (**)$$

Soit $\Pi(t) = \frac{1}{2i\pi}\int_{\Gamma_1}(z-Q(t))^{-1}\,dz$. Alors $\Pi(t)$ est un projecteur de rang 1 sur le sous-espace propre $\mathrm{Ker}(Q(t) - \lambda(t))$ tel que $\Pi(0)f = \nu(f)\mathbf{1}$, et l'on peut définir

$$v(t) = \langle \nu, \Pi(t)\mathbf{1}\rangle^{-1}\Pi(t)\mathbf{1} \quad \text{et} \quad \phi(t) = \Pi(t)^*\nu,$$

où $\Pi(t)^*$ est l'opérateur adjoint de $\Pi(t)$ (on verra ci-dessous que $\langle \nu, \Pi(t)\mathbf{1}\rangle \ne 0$ pour $|t|$ petit de sorte que $v(t)$ est bien défini). Le premier point de (b) est évident.

On a $\Pi(t)\mathbf{1} - \mathbf{1} = \Pi(t)\mathbf{1} - \Pi(0)\mathbf{1} = \frac{1}{2i\pi}\int_{\Gamma_1}[(z-Q(t))^{-1}\mathbf{1} - (z-Q)^{-1}\mathbf{1}]\,dz$, d'où

$$\nu(|\Pi(t)\mathbf{1} - \mathbf{1}|) \le \frac{1}{2\pi}\int_{\Gamma_1}\nu(|(z-Q(t))^{-1}\mathbf{1} - (z-Q)^{-1}\mathbf{1}|)\,dz \le C|t|^{\tau}.$$

On en déduit aisément (b1), ainsi que (b2) grâce à l'égalité $\langle \phi(t), \mathbf{1}\rangle = \langle \nu, \Pi(t)\mathbf{1}\rangle$. Enfin (b3) résulte de l'égalité

$$N(t)^n\mathbf{1} = N(t)^n\mathbf{1} - N(0)^n\mathbf{1} = \frac{1}{2i\pi}\int_{\Gamma_0} z^n[(z-Q(t))^{-1}\mathbf{1} - (z-Q)^{-1}\mathbf{1}]\,dz.$$

Il reste à prouver (c). Soit $t_0 \in J$. Les deux premières propriétés établies au début de la preuve et la propriété (H4) montrent que $Q(t)$ vérifie au voisinage de $t_0$ les conditions du théorème de Keller–Liverani. Ce dernier, avec la décomposition (D) écrite en $t = t_0$, assure que, pour $|h|$ petit, $Q(t_0 + h)$ admet une unique valeur propre dominante $z(t_0 + h)$ qui tend vers $\lambda(t_0)$ quand $h \to 0$. Mais par unicité on a $z(t_0 + h) = \lambda(t_0 + h)$. $\square$

En adaptant la preuve de [15], lemme 4.2, à l'aide de la majoration du corolaire 4.1, nous allons maintenant préciser le comportement de $\lambda(u)$ quand $u \to 0$.

**Lemme 4.2.** *Sous les hypothèses $(\mathcal{H})$ et $\sigma^2 > 0$, on a $\lambda(u) = 1 - \frac{\sigma^2}{2}u^2 + \mathrm{O}(|u|^{2+\tau})$ pour tout réel $0 < \tau < 1$. Si en outre $\langle \nu, |v(u) - \mathbf{1}|\rangle = \mathrm{O}(|u|)$, alors $\lambda(u) = 1 - \frac{\sigma^2}{2}u^2 + \mathrm{O}(u^3)$.*

**Démonstration.** On suppose pour simplifier que $\sigma^2 = 1$. En utilisant le fait que $(X_n)_n$ est une chaîne de Markov (voir par exemple [13], p. 23), puis le théorème 4.1, il est facile de voir que, pour tout $f \in \mathcal{B}$, pour toute probabilité initiale $\mu_0 \in \mathcal{B}'$, et pour $t \in J$, $n \ge 1$,

$$\mathbb{E}[f(X_n)\mathrm{e}^{\mathrm{i}tS_n}] = \langle \mu_0, Q(t)^n f\rangle = \lambda(t)^n\langle \phi(t), f\rangle\langle \mu_0, v(t)\rangle + \langle \mu_0, N(t)^n f\rangle. \qquad (***)$$

Dans la suite on considère $t \in [-1, 1]$ et un entier $n \ge N$, avec $N$ assez grand. Pour le moment $N$ est choisi tel que $\frac{t}{\sqrt{n}} \in J$ pour $n \ge N$. La formule $(***)$ appliquée en $\frac{t}{\sqrt{n}}$ avec $f = v(\frac{t}{\sqrt{n}})$ et $\mu_0 = \nu$ montre que

$$\lambda\left(\frac{t}{\sqrt{n}}\right)^n = \mathbb{E}\left[v\left(\frac{t}{\sqrt{n}}\right)(X_n)\mathrm{e}^{\mathrm{i}(t/\sqrt{n})S_n}\right].$$

---

[2] Les résultats de [19] sont présentés avec une norme auxiliaire $|\cdot|$ sur $\mathcal{B}$ vérifiant $|\cdot| \le \|\cdot\|$, mais on peut montrer que ceux-ci subsistent lorsque $|\cdot|$ est remplacée par une semi-norme, en l'occurrence ici $\nu(|\cdot|)$.

En utilisant l'inégalité triangulaire et l'invariance de $\nu$, il vient

$$\left|\lambda\left(\frac{t}{\sqrt{n}}\right)^n - e^{-t^2/2}\right| \leq \mathbb{E}\left[\left|v\left(\frac{t}{\sqrt{n}}\right)(X_n) - 1\right|\right] + \left|\mathbb{E}[e^{i(t/\sqrt{n})S_n}] - e^{-t^2/2}\right|$$

$$= \left\langle \nu, \left|v\left(\frac{t}{\sqrt{n}}\right) - \mathbf{1}\right|\right\rangle + |\mathbb{E}[e^{i(t/\sqrt{n})S_n}] - e^{-t^2/2}|.$$

Du corollaire 4.1 et du point (b1) du théorème 4.1 il vient pour $t \in [-1, 1]$ et $n \geq N$,

$$\left|\lambda\left(\frac{t}{\sqrt{n}}\right)^n - e^{-t^2/2}\right| \leq C\left|\frac{t}{\sqrt{n}}\right|^\tau.$$

Par ailleurs cette propriété est satisfaite avec $\tau = 1$ si $\langle \nu, |v(u) - \mathbf{1}|\rangle = O(|u|)$.

En utilisant la fonction log complexe définie pour $z \in \mathbb{C}$ non nul par $\log z = \ln|z| + i\arg(z)$, avec $\arg(z) \in\,]-\pi, \pi]$, on démontre qu'il existe une constante $C' > 0$ telle que l'on ait pour $t \in [-1, 1]$ et $n \geq N$, avec $N$ assez grand (voir les détails dans [15], p. 193, la continuité de $\lambda(\cdot)$ sur $J$ est importante pour ce point),

$$\left|n\log\lambda\left(\frac{t}{\sqrt{n}}\right) + \frac{t^2}{2}\right| \leq C'\left|\frac{t}{\sqrt{n}}\right|^\tau. \tag{$L$}$$

En observant maintenant que $z - 1 - \log z = (\log z)\alpha(z)$ avec $\alpha(z) = O(|z - 1|)$, on a

$$\left|n\left(\lambda\left(\frac{t}{\sqrt{n}}\right) - 1\right) + \frac{t^2}{2}\right| \leq n\left|\lambda\left(\frac{t}{\sqrt{n}}\right) - 1 - \log\lambda\left(\frac{t}{\sqrt{n}}\right)\right| + \left|n\log\lambda\left(\frac{t}{\sqrt{n}}\right) + \frac{t^2}{2}\right|$$

$$\leq \left|n\log\lambda\left(\frac{t}{\sqrt{n}}\right)\right|\left|\alpha\left(\lambda\left(\frac{t}{\sqrt{n}}\right)\right)\right| + C'\left|\frac{t}{\sqrt{n}}\right|^\tau.$$

En utilisant à nouveau ($L$) et le fait que $|t| \leq 1$, on voit que $|n\log\lambda(\frac{t}{\sqrt{n}})| \leq \frac{1}{2} + C'$.

En outre le corollaire 4.1 montre que la suite $(\frac{S_n}{\sqrt{n}})_n$ converge en loi vers $\mathcal{N}(0, \sigma^2)$. On déduit alors de [15], lemme 4.2, que $\lambda(u) = 1 - \frac{u^2}{2} + o(u^2)$. Donc $\alpha(\lambda(u)) = O(\lambda(u) - 1) = O(u^2)$. Par conséquent il existe une constante $C'' > 0$ telle que l'on ait pour $t \in [-1, 1]$ et $n \geq N$, avec $N$ assez grand,

$$\left|n\left(\lambda\left(\frac{t}{\sqrt{n}}\right) - 1\right) + \frac{t^2}{2}\right| \leq C''\left|\frac{t}{\sqrt{n}}\right|^\tau.$$

En divisant cette inégalité par $t^2$ pour $\frac{1}{2} \leq |t| \leq 1$, on obtient

$$\left|\left(\frac{t^2}{n}\right)^{-1}\left(\lambda\left(\frac{t}{\sqrt{n}}\right) - 1\right) + \frac{1}{2}\right| \leq \frac{C''}{t^2}\left|\frac{t}{\sqrt{n}}\right|^\tau \leq 4C''\left|\frac{t}{\sqrt{n}}\right|^\tau.$$

Soit $u \in \mathbb{R}^*$, $|u| \leq \frac{1}{\sqrt{N}}$. Il existe clairement $n \geq N$ tel que $\frac{1}{2\sqrt{n}} \leq |u| \leq \frac{1}{\sqrt{n}}$. L'inégalité précédente appliquée avec $t = \sqrt{n}u$ montre que $|\frac{\lambda(u) - 1}{u^2} + \frac{1}{2}| \leq 4C''|u|^\tau$. Ceci démontre le premier point du lemme.

Comme déjà indiqué, si $\langle \nu, |v(u) - \mathbf{1}|\rangle = O(|u|)$, les arguments précédents s'appliquent avec $\tau = 1$. $\qquad\square$

### 4.3. Démonstration du théorème 2.1

Supposons pour simplifier que $\sigma^2 = 1$, et rappelons que, en vertu de l'inégalité de Berry–Esseen, une vitesse en $n^{-\tau/2}$ ($\tau < 1$) sera obtenue dans le t.l.c. si l'on démontre que, pour un certain $\alpha > 0$, on a

$$A_n = \int_{-\alpha\sqrt{n}}^{\alpha\sqrt{n}}\left|\frac{\mathbb{E}[e^{it(S_n/\sqrt{n})}] - e^{-t^2/2}}{t}\right|dt = O(n^{-\tau/2}).$$

Pour le moment on choisit $\alpha > 0$ tel que $\alpha \in J$, où $J$ est l'intervalle du théorème 4.1.

En appliquant la formule $(***)$ de la section 4.2 avec $f = \mathbf{1}$, $\mu_0 = \nu$, et en posant $L(u) = \langle \phi(u), \mathbf{1} \rangle - 1$, on a $\mathbb{E}[\mathrm{e}^{\mathrm{i}t S_n / (\sqrt{n})}] = \lambda(\frac{t}{\sqrt{n}})^n + \lambda(\frac{t}{\sqrt{n}})^n L(\frac{t}{\sqrt{n}}) + \langle \nu, N(\frac{t}{\sqrt{n}})^n \mathbf{1} \rangle$. D'où

$$A_n \leq \int_{-\alpha\sqrt{n}}^{\alpha\sqrt{n}} \left| \frac{\lambda(t/\sqrt{n})^n - \mathrm{e}^{-t^2/2}}{t} \right| dt + \int_{-\alpha\sqrt{n}}^{\alpha\sqrt{n}} \left| \frac{\lambda(t/\sqrt{n})^n L(t/\sqrt{n})}{t} \right| dt + \int_{-\alpha\sqrt{n}}^{\alpha\sqrt{n}} \left| \frac{\langle \nu, N(t/\sqrt{n})^n \mathbf{1} \rangle}{t} \right| dt$$

$$= I_n + J_n + K_n.$$

Du lemme 4.2, on déduit que $|\lambda(u)| \leq 1 - \frac{u^2}{4} \leq \mathrm{e}^{-u^2/4}$ pour $|u|$ petit. Par conséquent on a pour $|\frac{t}{\sqrt{n}}| \leq \alpha$, avec $\alpha$ assez petit,

$$\left| \lambda\left(\frac{t}{\sqrt{n}}\right) \right| \leq \mathrm{e}^{-t^2/(4n)} \quad \text{d'où} \quad \left| \lambda\left(\frac{t}{\sqrt{n}}\right) \right|^n \leq \mathrm{e}^{-t^2/4}.$$

Comme dans la preuve du théorème de Berry–Esseen, on écrit

$$\lambda\left(\frac{t}{\sqrt{n}}\right)^n - \mathrm{e}^{-t^2/2} = \left( \lambda\left(\frac{t}{\sqrt{n}}\right) - \mathrm{e}^{-t^2/(2n)} \right) \sum_{k=0}^{n-1} \lambda\left(\frac{t}{\sqrt{n}}\right)^{n-k-1} \mathrm{e}^{-kt^2/(2n)}.$$

Il existe $C, C' > 0$ tels que l'ait pour $\frac{|t|}{\sqrt{n}}$ assez petit, d'une part $|\lambda(\frac{t}{\sqrt{n}}) - \mathrm{e}^{-t^2/(2n)}| \leq C |\frac{t}{\sqrt{n}}|^{2+\tau}$ (lemme 4.2), d'autre part $\sum_{k=0}^{n-1} |\lambda(\frac{t}{\sqrt{n}})|^{n-k-1} \mathrm{e}^{-kt^2/(2n)} \leq \sum_{k=0}^{n-1} \mathrm{e}^{-t^2(n-k-1)/(4n)-kt^2/(4n)} \leq C' n \mathrm{e}^{-t^2/4}$, d'où

$$\left| \lambda\left(\frac{t}{\sqrt{n}}\right)^n - \mathrm{e}^{-t^2/2} \right| \leq CC' n^{-\tau/2} |t|^{2+\tau} \mathrm{e}^{-t^2/4}.$$

La fonction $t \mapsto t^{1+\tau} \mathrm{e}^{-t^2/4}$ étant intégrable sur $\mathbb{R}$, on a $I_n = \mathrm{O}(n^{-\tau/2})$.

En outre, de $|L(\frac{t}{\sqrt{n}})| \leq C |\frac{t}{\sqrt{n}}|^\tau$ (théorème 4.1(b2)), et de l'intégrabilité de $t \mapsto |t|^{\tau-1} \mathrm{e}^{-t^2/4}$, on voit aisément que $J_n = \mathrm{O}(n^{-\tau/2})$. Enfin, du point (b3) du théorème 4.1, du fait que $u \mapsto u^{\tau-1}$ est intégrable sur $[-\alpha, \alpha]$ et que $\rho^n = \mathrm{O}(n^{-\tau/2})$, on obtient $K_n = \mathrm{O}(n^{-\tau/2})$.

## 4.4. *Démonstration du théorème 2.2*

On suppose ici que l'hypothèse $(\widetilde{\mathcal{H}})$ est satisfaite et que $\mu_0 \in \mathcal{B}' \cap \widetilde{\mathcal{B}}'$, où $\mu_0$ est la loi initiale. L'intégrale $A_n$ est définie comme en 4.3. En appliquant la formule $(***)$ avec $f = \mathbf{1}$, on obtient en posant ici $L(u) = \langle \phi(u), \mathbf{1} \rangle \langle \mu_0, v(u) \rangle - 1$,

$$A_n \leq \int_{-\alpha\sqrt{n}}^{\alpha\sqrt{n}} \left| \frac{\lambda(t/\sqrt{n})^n - \mathrm{e}^{-t^2/2}}{t} \right| dt + \int_{-\alpha\sqrt{n}}^{\alpha\sqrt{n}} \left| \frac{\lambda(t/\sqrt{n})^n L(t/\sqrt{n}))}{t} \right| dt + \int_{-\alpha\sqrt{n}}^{\alpha\sqrt{n}} \left| \frac{\langle \mu_0, N(t/\sqrt{n})^n \mathbf{1} \rangle}{t} \right| dt$$

$$= I_n + J_n + K_n.$$

Pour l'étude de $I_n$, $J_n$ et $K_n$, on considère des cercles orientés $\Gamma_1$ et $\Gamma_0$ comme dans la section 4.2. En utilisant la dernière condition de $(\widetilde{\mathcal{H}})$, il est facile de voir (cf. [15], p. 195) que, pour $f \in \mathcal{B}$, $t \in J$ et $z \in \Gamma_0 \cup \Gamma_1$, on a

$$\| (z - Q(t))^{-1} f - (z - Q)^{-1} f \|_\sim \leq C |t| \| f \|,$$

avec $C > 0$ indépendante de $f$, $t$ et $z$. Par intégration curviligne (voir les définitions de $v(t)$, $\phi(t)$, $N(t)^n$ dans la section 4.2), on obtient que

$\| v(u) - \mathbf{1} \|_\sim \leq C |u|$,
$| \langle \phi(u), \mathbf{1} \rangle - 1 | \leq C |u|$ (utiliser le fait que $\nu \in \widetilde{\mathcal{B}}'$),
$| \langle \mu_0, N(u)^n \mathbf{1} \rangle | \leq C \rho^n |u|$ (utiliser le fait que $\mu_0 \in \widetilde{\mathcal{B}}'$).

Le fait que $\widetilde{\mathcal{B}}$ s'envoie continûment dans $\mathbb{L}^1(\nu)$ montre que $\langle \nu, |v(u) - \mathbf{1}|\rangle = \mathrm{O}(|u|)$. Du lemme 4.2, il vient que $\lambda(u) = 1 - \frac{\sigma^2}{2}u^2 + \mathrm{O}(u^3)$. Les arguments vus à la section précédent pour la majoration de $I_n$ s'appliquent alors avec $\tau = 1$. Donc $I_n = \mathrm{O}(n^{-1/2})$. En outre, des majorations ci-dessus, on peut déduire que $L(u) = \mathrm{O}(|u|)$, et l'on a vu que $|\lambda(\frac{t}{\sqrt{n}})|^n \le \mathrm{e}^{-t^2/4}$ pour $|\frac{t}{\sqrt{n}}|$ assez petit ; on en déduit que $J_n = \mathrm{O}(n^{-1/2})$. Enfin on a clairement $K_n = \mathrm{O}(\rho^n) = \mathrm{O}(n^{-1/2})$.

## Appendice. Preuve de (H4) pour les chaînes géométriquement ergodiques

En conservant les hypothèses et notations du paragraphe 3.1, nous complétons ici la preuve du corollaire 3.1 en démontrant que $Q(t)$ vérifie (H4) sur l'espace $\mathcal{B}_W$.[3] Nous utiliserons la notation $r_{\mathrm{ess}}(\cdot)$ pour désigner le rayon spectral essentiel d'un opérateur.

La première inégalité de (H4) s'obtient facilement. En effet, pour $f \in \mathcal{B}_W$, on a $|Q(t)^n f| \le Q^n|f|$, donc $\|Q(t)^n f\|_W \le \|Q^n|f|\|_W$, et l'on conclut grâce à la propriété d'ergodicité $W$-géométrique.

Pour établir la seconde propriété de (H4), quitte à modifier le réel $\kappa$ dans (H4), il suffit de démontrer que $r_{\mathrm{ess}}(Q(t)) \le \rho < 1$, avec $\rho$ indépendant de $t$. À cet effet, on désignera par $\mathcal{B}^\infty$ l'espace des fonctions mesurables bornées sur $E$, muni de la norme de la convergence uniforme, et, étant donné un noyau $T(x, dy)$ sur $E$, on notera $\tilde{T}(x, dy) = W(x)^{-1}W(y)T(x, dy)$. Clairement, $T \in \mathcal{L}(\mathcal{B}_W)$ si et seulement si $\tilde{T} \in \mathcal{L}(\mathcal{B}^\infty)$, et $T$ et $\tilde{T}$ admettent le même spectre, en particulier le même rayon spectral essentiel, dans leur action respective sur $\mathcal{B}_W$ et $\mathcal{B}^\infty$. De l'ergodicité $W$-géométrique de $Q$, il vient : $\exists \ell \ge 1, \forall f \in \mathcal{B}_W, \|Q^\ell(f) - \nu(f)\mathbf{1}\|_W \le \frac{1}{2}\|f\|_W$.

Avec $f = 1_A W$ pour $A \in \mathcal{E}$ quelconque, on obtient $(\tilde{Q})^\ell(x, A) = \widetilde{Q^\ell}(x, A) \le \frac{1}{2} + \nu(1_A W)$ pour tout $x \in E$. Définissons maintenant la probabilité $\widetilde{\nu}(A) = \nu(W)^{-1}\nu(W1_A)$. On a

$$\forall A \in \mathcal{E}, \quad \left(\widetilde{\nu}(A) \le \frac{1}{4\nu(W)}\right) \quad \Longrightarrow \quad \left(\forall x \in E, (\tilde{Q})^\ell(x, A) \le \frac{3}{4}\right).$$

Soit $t \in \mathbb{R}$. Avec les notations de [12], posons $\chi(y) = \mathrm{e}^{\mathrm{i}t\xi(y)}$ et $\tilde{Q}_\chi(x, dy) = \chi(y)\tilde{Q}(x, dy)$. En appliquant [12] (lemma 3.4, theorem 3.1), l'implication ci-dessus montre que $r_{\mathrm{ess}}(\tilde{Q}_\chi) \le (\frac{3}{4})^{1/\ell}$. Puisque $\tilde{Q}_\chi = \widetilde{Q(t)}$, on obtient $r_{\mathrm{ess}}(Q(t)) = r_{\mathrm{ess}}(\widetilde{Q(t)}) \le (\frac{3}{4})^{1/\ell}$.

**Nota.** *Actuellement au moment de la parution de cet article, une vitesse en $n^{-1/2}$ dans le t.l.c. a été obtenue, sous la condition optimale $\nu(|\xi|^3) < +\infty$, pour les chaînes uniformément ergodiques (i.e. cas $V = 1_E$ de la section 3.1) qui correspondent aux chaînes apériodiques vérifiant la condition de Doeblin. Plus précisément, pour une telle chaîne, la condition $(\widetilde{\mathcal{H}})$ peut être établie avec $\mathcal{B} = \mathbb{L}^3(\nu)$ et $\widetilde{\mathcal{B}} = \mathbb{L}^{3/2}(\nu)$ ; le théorème 2.2 s'applique alors lorsque $\xi \in \mathbb{L}^3(\nu)$. Cet exemple sera détaillé dans un prochain travail, fait en collaboration avec F. Pène, où nous établissons en outre, dans le cadre général des chaînes fortement ergodiques, un développement d'Edgeworth d'ordre 1 et un résultat de vitesse de convergence dans le t.l.c. multi-dimensionnel.*

---

[3]Mentionnons que cette propriété permet d'étendre le théorème local de [15] (Prop. 3.1) à toute chaîne géométriquement ergodique (la conditon restrictive $Q(x, dy) = K(x, y) \, d\lambda(y)$ dans [15] n'est pas nécessaire).